\documentclass[12pt]{amsart}
\usepackage[mathscr]{eucal}
\usepackage{amssymb}
\usepackage{latexsym}
\usepackage{amsthm}
\usepackage{amsmath}
\usepackage{graphicx}

\theoremstyle{plain}

\newtheorem{thm}{Theorem}[section]
\newtheorem{cor}{Corollary}[section]
\newtheorem{lem}{Lemma}[section]

\theoremstyle{definition}

\newtheorem{prf}{Proof}[section]

\title{Summability methods and Fourier analysis on $\mathbb{R}^{\times}$}
\author{Ryoichi Kunisada}

\address{Faculty of Education and Integrated Arts and Science, Waseda University, Shinjuku-ku, Tokyo 169-8050, Japan}

\email{rkunisada@aoni.waseda.jp}

\date{}

\begin{document}
\maketitle

\begin{abstract}
We introduce a certain class of summability methods on $L^{\infty}([1, \infty))$ which is defined by convolution in the group algebra $L^1(\mathbb{R}^{\times})$. This class contains an integral version of Ces\`{a}ro summability method and we particularly give a necessary and sufficient condition for a summability method in the class to equivalent to this one. Also, we consider another class of summability methods concerned with convolution in the group algebra $L^1(\mathbb{R})$ and give similar results.
\end{abstract}

\bigskip

\section{Introduction}
Let $\mathbb{R}_+ = [0, \infty)$. Let us $L^{\infty}(\mathbb{R}_+)$ be the set of all essentially bounded measurable functions on $\mathbb{R}_+$ and $L^1(\mathbb{R}_+)$ be the set of all integrable functions on $\mathbb{R}_+$ (with respect to the Lebesgue measure dt). Recall that group algebra $L^1(\mathbb{R})$ of $\mathbb{R}$ is a Banach algebra when multiplication is defined by convolution
\[(f*g)(x) = \int_{-\infty}^{\infty} f(x-t)g(t)dt, \]
where $f, g \in L^1(\mathbb{R})$.

Throughout the paper, by a summability method on $L^{\infty}(\mathbb{R}_+)$ we mean a continuous linear functional $S$ defined on a certain closed subspace $\mathcal{D}(S)$ of $L^{\infty}(\mathbb{R}_+)$. When we introduce a summability method $(S, \mathcal{D}(S))$, we will refer to only the expression of $S$, which involves a certain limit process, and let us its domain $\mathcal{D}(S)$ be the closed subspace consisting of all $f \in L^{\infty}(\mathbb{R}_+)$ for which the limit in the definition of $S$ exists. 
For given summability methods $(S, \mathcal{D}(S))$ and $(S^{\prime}, \mathcal{D}(S^{\prime}))$, we say that $S$ is stronger than $S^{\prime}$ if both $\mathcal{D}(S^{\prime}) \subseteq \mathcal{D}(S)$ and $S=S^{\prime}$ on $\mathcal{D}(S^{\prime})$ holds. Further, we say that $(S, \mathcal{D}(S))$ and $(S^{\prime}, \mathcal{D}(S^{\prime}))$ are equivalent if the converse is also true. We also assume that the summability methods are regular, namely, if $\lim_{x \to \infty} f(x) = \alpha$ then $S(f) = \alpha$ holds.

First we introduce the summability method $K$ on $L^{\infty}(\mathbb{R}_+)$ which is defined by
\[K(f) = w^*\mathchar`-\lim_s f(x+s), \]
where the limit in the definition means the weak* convergence in $L^{\infty}(\mathbb{R}_+)$. In other words, for any $f(x) \in L^{\infty}(\mathbb{R}_+)$, $K(f) = \alpha$ if and only if $f(x+s)$ converges to the constant function $\alpha$ as $s$ tends to $\infty$ with respect to the weak* topology of $L^{\infty}(\mathbb{R}_+)$.
In particular, one of our main objectives is the summability methods of the following type. For each $\phi \in L^1(\mathbb{R}_+)$ with $\int_0^{\infty} \phi(x)dx = 1$, which is imposed to guarantee that $S_{\phi}$ is regular, we define a summability method $S_{\phi}$ as
\[S_{\phi}(f) = \lim_{x \to \infty} \int_0^x f(t)\phi(x-t)dt, \]
where $f \in L^{\infty}(\mathbb{R}_+)$. We denote the set of all such summability methods by $\Phi$. One of our main aim is to give a necessary and sufficient condition for a given $S_{\phi} \in \Phi$ to equivalent to $K$. Moreover, we will give a sufficient condition that for any pair of summability methods in $\Phi$ one is stronger than the other one.

Now we mention the relation of $\Phi$ with Wiener's Tauberian theorem. Let $\tilde{f}$ and $\tilde{\phi}$ be the functions on $\mathbb{R}$ which extend $f$ and $\phi$ respectively by defining $\tilde{f} = \tilde{\phi} = 0$ on the interval $(-\infty, 0)$. Then notice that it can be written as
\[S_{\phi}(f) = \lim_{x \to \infty} \int_0^x f(t)\phi(x-t)dt = \lim_{x \to \infty} \int_{-\infty}^{\infty} \tilde{f}(t)\tilde{\phi}(x-t)dt = \lim_{x \to \infty} (\tilde{f}*\tilde{\phi})(x) \]
since the integrand vanishes on the intervals $(-\infty, 0)$ and $(x, \infty)$. In this sense, the summability methods in $\Phi$ can be viewed as the special case of Wiener's Tauberian theorem in which both $f(x) \in L^{\infty}(\mathbb{R})$ and $\phi(x) \in L^1(\mathbb{R})$ vanishes on the half line $(-\infty, 0)$.

We also consider the multiplicative version, in which we are more interested in this paper. Remark that in the following definitions, the underlying measure is the Haar measure of the multiplicative group $\mathbb{R}^{\times} = (0, \infty)$ of $\mathbb{R}$, namely, $dt/t$.
Let $\mathbb{R}_+^{\times} = [1, \infty)$. Let us $L^{\infty}(\mathbb{R}_+^{\times})$ be the set of all essentially bounded measurable functions on $\mathbb{R}_+^{\times}$ and $L^1(\mathbb{R}_+^{\times})$ be the set of all integrable functions on $\mathbb{R}_+^{\times}$ (with respect to $dt/t$). Notice that the group algebra $L^1(\mathbb{R}^{\times})$ of $\mathbb{R}^{\times}$ is a Banach algebra with the multiplication
\[(f*g)(x) = \int_0^{\infty} f(x/ t)g(t)\frac{dt}{t},  \]
where $f, g \in L^1(\mathbb{R}^{\times})$. 

The definition of summability methods on $L^{\infty}(\mathbb{R}_+^{\times})$ is obtained by simply replacing the word `$L^{\infty}(\mathbb{R}_+)$' by `$L^{\infty}(\mathbb{R}_+^{\times})$' in the definition of summability methods on $L^{\infty}(\mathbb{R}_+)$. Other definitions concerned with summability methods on $L^{\infty}(\mathbb{R}_+)$ can also be transferred to this context.
Similarly, we introduce the summability method $P$ on $L^{\infty}(\mathbb{R}_+^{\times})$ by
\[P(f) = w^*\mathchar`-\lim_r f(rx), \]
where the limit in the definition means the weak* convergence in $L^{\infty}(\mathbb{R}_+^{\times})$.
Also for any $\psi \in L^1(\mathbb{R}_+^{\times})$ with $\int_1^{\infty} \psi(x)\frac{dx}{x} = 1$, we consdier a summability method $S_{\psi}$ defined as
\[M_{\psi}(f) = \lim_{x \to \infty} \int_1^x  f(t)\psi(x/t) \frac{dt}{t}, \]
where $f \in L^{\infty}(\mathbb{R}_+^{\times})$. We denote the set of all such summability methods by $\Psi$. Then we will also obtain a necessary and sufficient condition for a given $M_{\psi}$ to equivalent to $P$.
In particular, we consider the summability method $M$ defined by
\[M(f) = \lim_{x \to \infty} \frac{1}{x} \int_1^x f(t)dt, \]
which is a continuous version of Ces\`{a}ro summability method. Notice that it is an element of $\Psi$. In fact, it can be written as
\[M(f) = \lim_{x \to \infty} \int_1^x f(t) \left(\frac{x}{t}\right)^{-1} \frac{dt}{t}, \]
where $\psi(x) = x^{-1}$ is in $L^1(\mathbb{R}_+^{\times})$. We show that $M$ is equivalent to $P$ as an application of the above result, which means that we will also obtain a necessary and sufficient condition that $M_{\psi} \in \Psi$ is equivalent to $M$.

\section{Preliminaries}
The proof of  our main theorems are based on the notion of the limit along an ultrafilter, whose definition is as follows: let $f : X \rightarrow Y$ be a mapping of a set $X$ into a compact space $Y$ and $\mathcal{U}$ be an ultrafilter on $X$. Then there exists an element $y$ of $Y$ such that $f^{-1}(U) \in \mathcal{U}$ holds for every neighborhood $U$ of $y$. This element $y$ of $Y$ is called the limit of $f$ along $\mathcal{U}$ and denoted by $\mathcal{U}\mathchar`-\lim_x f(x)$. In this paper we will consider only the case $X = \mathbb{N}$, $\mathbb{R}_+$, or $\mathbb{R}_+^{\times}$. In particular, we consider an ultrafilter on $\mathbb{N}$ as an element of the Stone-\v{C}ech compactification $\beta\mathbb{N}$ of $\mathbb{N}$. Recall that $\beta\mathbb{N}$ is a compactification of $\mathbb{N}$ with the property that each mapping $f : \mathbb{N} \rightarrow X$ of $\mathbb{N}$ into a compact space $X$ can be extended to a continuous mapping $\overline{f} : \beta\mathbb{N} \rightarrow X$ of $\beta\mathbb{N}$ into $X$. In fact, for each $\eta \in \beta\mathbb{N}$, $\overline{f}(\eta)$ is given by the formula $\overline{f}(\eta) = \eta\mathchar`-\lim_n f(n)$. In particular, for any $f \in l_{\infty}$ of the set of all bounded functions on $\mathbb{N}$, its continuous extension to $\beta\mathbb{N}$ gives an algebraic isomorphism between $l_{\infty}$ and $C(\beta\mathbb{N})$ of the set of all continuous functions on $\beta\mathbb{N}$, which implies that $\beta\mathbb{N}$ is the maximal ideal space of $l_{\infty}$.

Of particular importance concerning $\beta\mathbb{N}$ is the extension of the translation of $\mathbb{N}$; let $\tau_0 : \mathbb{N} \rightarrow \mathbb{N}$ be a mapping $n \mapsto n+1$, Viewing $\tau_0$ as a mapping of $\mathbb{N}$ into $\beta\mathbb{N}$, one can extend it continuously to $\tau : \beta\mathbb{N} \rightarrow \beta\mathbb{N}$ of a continuous mapping of $\beta\mathbb{N}$ into itself. In particular, let us $\mathbb{N}^* = \beta\mathbb{N} \setminus \mathbb{N}$ be the set of all free-ultrafilters on $\mathbb{N}$ and restrict $\tau$ to $\mathbb{N}^*$, we get a homeomorphism $\tau : \mathbb{N}^* \rightarrow \mathbb{N}^*$ of $\mathbb{N}^*$ onto itself. Namely, the pair $(\mathbb{N}^*, \tau)$ is a discrete flow. 

Next we introduce continuous counterparts of these notions. Let $C_{ub}(\mathbb{R}_+)$ be the set of all uniformly continuous bounded functions on $\mathbb{R}_+$. The maximal ideal space $\Omega$ of $C_{ub}(\mathbb{R}_+)$ is a compact space given by the following construction (see [4] for a proof): let us consider the product space $\beta\mathbb{N} \times [0, 1]$ and define $\Omega$ as a quotient space by identifying the points $(\eta, 1)$ and $(\tau\eta, 0)$ in $\beta\mathbb{N} \times [0, 1]$ for every $\eta \in \beta\mathbb{N}$. Notice that $\Omega$ is a compactification of $\mathbb{R}_+$ to which every uniformly continuous bounded function on $\mathbb{R}_+$ can be extended continuously. For each $f \in C_{ub}(\mathbb{R}_+)$ and $\omega = (\eta, u) \in \Omega$, the extended function $\overline{f}(\omega) \in C(\Omega)$ is given by
\[\overline{f}(\omega) = \omega\mathchar`-\lim_x f(x), \]
where we identify $\omega$ with the ultrafilter on $\mathbb{R}_+$ defined by $\{u+A : A \in \eta\}$.
We also consider a continuous version of the discrete flow $(\mathbb{N}^*, \tau)$. Let us denote by $\Omega^*$ the closed subspace of $\Omega$ consisting of all elements $\omega = (\eta, u)$ in $\Omega$ with $\eta  \in \mathbb{N}^*$. Then we define a continuous flow on $\Omega^*$ as follows; for each $s \in \mathbb{R}$, we define the homeomorphism $\tau^s : \Omega^* \rightarrow \Omega^*$ by
\[\tau^s(\eta, u) = (\tau^{[u+s]}\eta, u+s-[u+s]), \]
where $[x]$ denotes the largest integer not exceeding $x$ for a real number $x$. Notice that the continuous flow $(\Omega^*, \{\tau^s\}_{s \in \mathbb{R}})$ is the suspension of the discrete flow $(\mathbb{N}^*, \tau)$. 

For each $\omega \in \Omega^*$, we denote by $o(\omega) = \{\tau^x \omega : x \in \mathbb{R}\}$ the orbit of $\omega$. For any $f \in C_{ub}(\mathbb{R}_+)$, let $f_{\omega}(x) = \overline{f}(\tau^x \omega)$ be the restriction of $\overline{f}(\omega)$ to the orbit $o(\omega)$ of $\omega \in \Omega$, which is a uniformly continuous bounded function on $\mathbb{R}$. Now we extend this mapping $C_{ub}(\mathbb{R}_+) \ni f(x) \mapsto f_{\omega}(x) \in C_{ub}(\mathbb{R})$ to $L^{\infty}(\mathbb{R}_+)$. We write $f_s(x) = f(x+s)$ for every $f(x) \in L^{\infty}(\mathbb{R}_+)$ and $s \ge 0$. For any given $f(x) \in L^{\infty}(\mathbb{R}_+)$, since $\{f_s(x)\}_{s \ge 0}$ is a uniformly bounded family of $L^{\infty}(\mathbb{R}_+)$, it is a weak* relatively compact subset of $L^{\infty}(\mathbb{R}_+)$ and hence we can consider its limit along an ultrafilter $\omega \in \Omega^*$. We define $f_{\omega}(x) = \omega\mathchar`-\lim_s f_s(x)$ for each $\omega \in \Omega^*$. Now we can extend it to the negative direction in a natural way as
\[f_{\omega}(x) = f_{\tau^{-N}\omega}(N+x), \quad x \in [-N, 0], \]
for every $N > 0$. In this way, we consider $f_{\omega}$ to be a function in $L^{\infty}(\mathbb{R})$.
In particular, it is easy to see that $f_{\omega}(x) = \overline{f}(\tau^x \omega)$ holds if $f(x)$ is in $C_{ub}(\mathbb{R}_+)$. 

Now the following lemma is obvious.

\begin{lem}
For any $f(x) \in C_{ub}(\mathbb{R}_+)$, $\lim_{x \to \infty} f(x) = \alpha$ if and only if $\overline{f}(\omega) = \alpha$ for every $\omega \in \Omega^*$.
\end{lem}

Now we partition $\Omega^*$ into the disjoint orbits and denote by $\mathscr{O}$ the set of all representatives of each orbit in $\Omega^*$. Then Lemma 2.1 can be restated as follows.

\begin{lem}
For any $f(x) \in C_{ub}(\mathbb{R}_+)$, $\lim_{x \to \infty} f(x) = \alpha$ if and only if $f_{\omega}(x) = \alpha$ for every $\omega \in \mathscr{O}$.
\end{lem}

We can generalize the above lemma to the case of $f(x) \in L^{\infty}(\mathbb{R}_+)$.
\begin{lem}
For any $f(x) \in L^{\infty}(\mathbb{R}_+)$, $K(f) = \alpha$ if and only if $f_{\omega}(x) = \alpha$ for every $\omega \in \mathscr{O}$.
\end{lem}

\begin{prf}
The necessity is obvious. Suppose that $f_{\omega}(x) = \alpha$ for every $\omega \in \mathcal{O}$. Then for any $\phi \in L^1(\mathbb{R}_+)$ and $\omega \in \Omega^*$, we have
\[\omega\mathchar`-\lim_x \int_0^{\infty} f_x(t)\phi(t)dt = \int_0^{\infty} f_{\omega}(t)\phi(t)dt = \alpha \int_0^{\infty} \phi(t)dt. \]
Hence by Lemma 2.1, we have
\[\lim_{x \to \infty} \int_0^{\infty} f_x(t)\phi(t)dt = \alpha \int_0^{\infty} \phi(t)dt \]
for every $\phi \in L^1(\mathbb{R}_+)$, which means that $w^*\mathchar`-\lim_x f_x(t) = \alpha$, that is, $K(f) = \alpha$.
\end{prf}

Similarly, for any given $f(x) \in L^{\infty}(\mathbb{R}_+^{\times})$, let us consider the uniformly bounded family $\{f(rx)\}_{r \ge 1}$ of $L^{\infty}(\mathbb{R}_+^{\times})$ and define $f_{\omega}^{\times}(x) = \omega\mathchar`-\lim_s f(e^sx)$ for each $\omega \in \Omega^*$, where the limit is taken with respect to the weak* topology of $L^{\infty}(\mathbb{R}_+^{\times})$. Notice that it can be written also as $f_{\omega}^{\times}(x) = e^{\omega}\mathchar`-\lim_r f(rx)$, where $e^{\omega}$ is an ultrafilter on $\mathbb{R}_+^{\times}$ consisting of all sets $e^X = \{e^x : x \in X\}$ for every $X \in \omega$. Also we consider $f_{\omega}^{\times}(x)$ to be a function in $L^{\infty}(\mathbb{R}^{\times})$ in a similar way as $f_{\omega}(x)$. Then we have an analogous result of Lemma 2.3.

\begin{lem}
For any $f(x) \in L^{\infty}(\mathbb{R}_+^{\times})$, $P(f) = \alpha$ if and only if $f_{\omega}^{\times}(x) = \alpha$ for every $\omega \in \mathscr{O}$.
\end{lem}

We introduce some notions and results from Fourier analysis on $\mathbb{R}$. Recall that for a function $\phi \in L^1(\mathbb{R})$, its Fourier transform $\hat{\phi} \in C_0(\mathbb{R})$ is defined by
\[\hat{\phi}(x) = \int_{-\infty}^{\infty} \phi(t)e^{-ixt}dt, \quad x \in \mathbb{R}. \]
For a function $\phi \in L^1(\mathbb{R})$, let us $Z(\phi) = \{x \in \mathbb{R} : \hat{\phi}(x) = 0\}$. Suppose that $I$ be a closed ideal of $L^1(\mathbb{R})$. Then we define the zero-set of $I$ by $Z(I) = \cap_{\phi \in I} Z(\phi)$. Notice that $L^1(\mathbb{R})^* = L^{\infty}(\mathbb{R})$, where we identify each $f \in L^{\infty}(\mathbb{R})$ as the linear functional on $L^1(\mathbb{R})$ defined  by
\[L^1(\mathbb{R}) \ni \phi(x) \mapsto \int_{-\infty}^{\infty} \phi(x)f(-x)dx. \]
Let us $A = I^{\bot}$, the annihilator of $I$, that is, $A = \{f(x) \in L^{\infty}(\mathbb{R}) : \int_{-\infty}^{\infty} f(x)\phi(-x)dx = 0, \ \forall \phi \in I \}$. Then $A$ is a weak* closed translation invariant subspace of $L^{\infty}(\mathbb{R})$. The spectrum of $A$ is defined by $\sigma(A) = \{\alpha \in \mathbb{R} : e^{i\alpha x} \in A\}$. In particular, it holds that $\sigma(A) = Z(I)$ (See [6] for more details of these notions). 

Remark that the notions of the zero-set of a closed ideal in $L^1(\mathbb{R})$ and the spectrum of a weak* closed translation invariant subspace in $L^{\infty}(\mathbb{R})$ can be defined to the multiplicative group $\mathbb{R}^{\times}$ of $\mathbb{R}$ or more generally, to arbitrary locally compact abelian groups. See for example [6].

\section{Main results}
First we work in the additive setting. For a given $\phi(x) \in L^1(\mathbb{R}_+)$, let us define the linear operator $U_{\phi}$ by
\[U_{\phi} : L^{\infty}(\mathbb{R}_+) \longrightarrow L^{\infty}(\mathbb{R}_+), \quad (U_{\phi}f)(x) = \int_0^x f(t)\phi(x-t)dt. \]

\begin{lem}
For any $f \in L^{\infty}(\mathbb{R}_+)$ and $\phi \in L^1(\mathbb{R}_+)$, $(U_{\phi}f)(x)$ is in $C_{ub}(\mathbb{R}_+)$.
\end{lem}

\begin{prf}
For any $x, y \in \mathbb{R}_+$ with $x < y$, we have
\begin{align}
|(Uf)(x) - (Uf)(y)| &= \left|\int_0^x f(t)\phi(x-t)dt - \int_0^y f(t)\phi(y-t)dt \right| \notag \\
&\le \left|\int_0^x f(t)(\phi(x-t)-\phi(y-t))dt \right| + \int_x^y |f(t)| \cdot |\phi(y-t)|dt \notag \\
&\le \|f\|_{\infty} \int_0^x |\phi(x-t)-\phi(y-t)|dt + \|f\|_{\infty} \int_x^y |f(t)| \cdot |\phi(y-t)|dt \notag \\
&\le \|f\|_{\infty} \int_0^{\infty} |\phi(t) - \phi(t+y-x)|dt + \|f\|_{\infty} \int_0^{y-x} |\phi(t)|dt. \notag
\end{align}
Therefore we get the result since the last two definite integrals tend to $0$ as $|x-y| \rightarrow 0$.
\end{prf}

Then the next result plays an essential role for the remainder of this paper.

\begin{lem}
For every $\phi \in L^1(\mathbb{R}_+)$, $f \in L^{\infty}(\mathbb{R}_+)$ and $\omega \in \Omega^*$, it holds that
\[\omega\mathchar`-\lim_x \int_0^x f(t)\phi(x-t)dt = \int_0^{\infty} f_{\omega}(-t)\phi(t)dt. \]
\end{lem}

\begin{prf}
We have
\begin{align}
\omega\mathchar`-\lim_x \int_0^x f(t)\phi(x-t)dt &= \omega\mathchar`-\lim_x \int_0^x f(x-t)\phi(t)dt = 
\omega\mathchar`-\lim_x \int_0^x f_x(-t)\phi(t)dt \notag \\
&= \omega\mathchar`-\lim_x \int_{-\infty}^{\infty} \tilde{f_x}(-t)\tilde{\phi}(t)dt = \int_{-\infty}^{\infty} \tilde{f_{\omega}}(-t)\tilde{\phi}(t)dt  \notag \\
&= \int_0^{\infty} f_{\omega}(-t)\phi(t)dt, \notag
\end{align}
where $\tilde{f}_{\omega}(-t)$ is the weak* limit of $\{\tilde{f}_x(-t)\}_{x \ge 0}$ along $\omega$ in the space $L^{\infty}(\mathbb{R})$. The proof is complete.
\end{prf}

\begin{thm}
For every $\phi \in L^1(\mathbb{R}_+)$, $f \in L^{\infty}(\mathbb{R}_+)$ and $\omega \in \Omega^*$, we have
\[(U_{\phi}f)_{\omega}(x) = (f_{\omega} * \tilde{\phi})(x) \]
for every $x \in \mathbb{R}$, where $\tilde{\phi}$ is a function in $L^1(\mathbb{R})$ such that
\[\tilde{\phi}(x) = \begin{cases}
            \phi(x) & \text{if $x \ge 0$}, \\
            0       & \text{if $x < 0$}.
         \end{cases}     
\]
\end{thm}

\begin{prf}
Observe that
\begin{align} 
(U_{\phi}f)_{\omega}(x) &= \tau^x\omega\mathchar`-\lim_s \int_0^s f(t)\phi(s-t)dt = \int_0^{\infty} f_{\tau^x \omega}(-t) \phi(t)dt  \notag \\
&= \int_0^{\infty} f_{\omega}(x-t)\phi(t)dt = (f_{\omega} * \tilde{\phi})(x). \notag 
\end{align}
\end{prf}

\begin{thm}
For every $\phi \in L^1(\mathbb{R}_+)$ with $\int_0^{\infty}\phi(x)dx = 1$, the summability method $S_{\phi}$ is equivalent to the summability method $K$ if and only if the closed ideal $I$ generated by $\tilde{\phi}$ in $L^1(\mathbb{R})$ is the whole space $L^1(\mathbb{R})$.
\end{thm}

\begin{prf}
By the assumption that $\int_0^{\infty} \phi(t)dt = 1$, for any $\omega \in \Omega^*$ we have
\begin{align}
(f_{\omega} * \tilde{\phi})(x) = \alpha &\Longleftrightarrow (f_{\omega} - \alpha) * \tilde{\phi} = 0 \notag \\
&\Longleftrightarrow \int_{-\infty}^{\infty} (f_{\omega} - \alpha)(t)\tilde{\phi}(x-t)dt = 0 \notag \\
&\Longleftrightarrow \int_{-\infty}^{\infty} (f_{\omega} - \alpha)(t) \tilde{\phi}_x(-t)dt = 0 \notag \\
&\Longleftrightarrow \int_{-\infty}^{\infty} (f_{\omega} - \alpha)(t) \phi_1(-t)dt = 0, \notag
\end{align}
where $\phi_1$ is any element of the closed subspace generated by all translates $\tilde{\phi}_x$ of $\tilde{\phi}$, which is equal to the closed ideal $I$ of $L^1(\mathbb{R})$ generated by $\tilde{\phi}$. Hence if $I = L^1(\mathbb{R})$, then
by the Hahn-Banach theorem, $f_{\omega} - \alpha = 0$ holds for every $\omega \in \Omega^*$. This shows that $K(f) = \alpha$ by Lemma 2.3. On the other hand, if $I \subsetneq L^1(\mathbb{R})$, then $Z(I) \not= \emptyset$. Let us $A = I^{\bot}$ and we have again by the Hahn-Banach theorem, that $A \not= \{0\}$. In particular, we have $\sigma(A) \not= \emptyset$. Thus there exists some $\alpha \in \mathbb{R}$ such that $e^{i\alpha x} \in A$. Put $f(x) = e^{i \alpha x}$ and notice that for every $\omega \in \Omega^*$, it holds that $f_{\omega}(x) = f_{\omega}(0)e^{i \alpha x}$. By the assumption that $A = I^{\bot}$, we have
\[\overline{(U_{\phi}f)}(\omega) = \int_{-\infty}^{\infty} f_{\omega}(-t)\tilde{\phi}(t)dt = \int_{-\infty}^{\infty} f_{\omega}(0)e^{-i \alpha t}\tilde{\phi}(t)dt = f_{\omega}(0)\hat{\tilde{\phi}}(\alpha) = 0 \]
for every $\omega \in \Omega^*$, which indicates $S_{\phi}(f) = 0$. But $f_{\omega}(x) = f_{\omega}(0)e^{i \alpha x} \not= 0$ for every $\omega \in \Omega^*$, i.e., $K(f) \not= 0$. Thus $S_{\phi}$ is not equivalent to $K$.
\end{prf}

Combining Theorem 3.2 with Wiener's Tauberian theorem, we get the following result immediately.
\begin{cor}
For every $\phi \in L^1(\mathbb{R}_+)$ with $\int_0^{\infty} \phi(x)dx = 1$, the summability method $S_{\phi}$ is equivalent to the summability method $K$ if and only if $Z(\tilde{\phi}) = \emptyset$.
\end{cor}

The following corollary is obvious by the proof of Theorem 3.2.
\begin{cor}
$K$ is the weakest summability method in $\Phi$.
\end{cor}

In particular, Theorem 3.2 can be viewed as the comparison theorem between the weakest element $K$ in $\Phi$ and any $S_{\phi}$ in $\Phi$. For a general pair $(S_{\phi}, \mathcal{D}(S_{\phi}))$ and $(S_{\phi^{\prime}}, \mathcal{D}(S_{\phi^{\prime}}))$ of $\Phi$, the following result is immediate by the proof of the sufficiency of Theorem 3.2.

\begin{thm}
Let $I_{\phi}$ and $I_{\phi^{\prime}}$ be the closed ideals in $L^1(\mathbb{R})$ generated by $\tilde{\phi}$ and $\tilde{\phi^{\prime}}$ respectively. Then $S_{\phi}$ is stronger than $S_{\phi^{\prime}}$ if $I_{\phi} \subseteq I_{\phi^{\prime}}$.
\end{thm}

Identifying $\phi(x)$ in $L^1(\mathbb{R}_+)$ with $\tilde{\phi}(x)$ in $L^1(\mathbb{R})$, we can view $L^1(\mathbb{R}_+)$ as a Banach algebra since for any pair $\phi_1, \phi_2 \in L^1(\mathbb{R}_+)$ the convolution $\tilde{\phi}_1 * \tilde{\phi}_2$ also vanishes on the half line $(-\infty, 0)$ and can be written by $\tilde{\phi}_3$ for some $\phi_3 \in L^1(\mathbb{R}_+)$. From now on, we wil simply denote it as $\phi_1 * \phi_2 = \phi_3$.

\begin{thm}
Let $\phi_1, \phi_2 \in L^1(\mathbb{R}_+)$. Consider the summability method $S_{\phi_2 \circ \phi_1}$ defined by
\[S_{\phi_2 \circ \phi_1} = \lim_{x \to \infty} (U_{\phi_2}U_{\phi_1}f)(x). \]
Then $S_{\phi_2 \circ \phi_1} = S_{\phi_2 * \phi_1}$ holds.
\end{thm}

\begin{prf}
By associativity of the convolution product and the fact that $(U_{\phi}f)(x) = (\tilde{f}*\tilde{\phi})(x)$ as we have mentioned in Section 1, the assertion follows immediately.
\end{prf}

The following result follows immediately by Theorem 3.4.
\begin{cor}
Let $\phi_1, \phi_2, \ldots, \phi_k$ be a set of functions in $L^1(\mathbb{R}_+)$ with $\int_0^{\infty}\phi_i(x)dx = 1$ for every $1 \le i \le k$ such that the closed ideal generated by $\tilde{\phi}_1, \tilde{\phi}_2, \ldots, \tilde{\phi}_k$ is $L^1(\mathbb{R})$. Then the assertion that $S_{\phi_i}(f) = \alpha$ for every $1 \le i \le k$ is equivalent to the assertion that $K(f) = \alpha$.
\end{cor}

In particular, we consider iterations of a summability method $S_{\phi}$, namely, for every $k \ge 1$, we define the summability methods $S_{\phi}^k$ by
\[S_{\phi}^k(f) = \lim_{x \to \infty} (U_{\phi}^kf)(x). \]

\begin{cor}
If $S_{\phi}$ is equivalent to $K$, then every $S_{\phi}^k$ is also equivalent to $K$, that is, $S_{\phi}$ and $S_{\phi}^k$ is equivalent for every $k \ge 1$.
\end{cor}

\begin{prf}
By Theorem 3.4, for every $k \ge 1$, we have $S_{\phi}^k = S_{{\phi}^k}$.
Notice that the Fourier transform of the convolution product $\phi_1 * \phi_2$ of $\phi_1$ and $\phi_2$ in $L^1(\mathbb{R})$ is the pointwise product $\hat{\phi}_1 \cdot \hat{\phi}_2$ of their Fourier transforms. Hence $\hat{\phi^k} = (\hat{\phi})^k$ and it means that if $\hat{\phi}$ has no zero point, then $\hat{\phi^k}$ has also no zero point. Thus by Corollary 3.1, we get the result.
\end{prf}

Next we deal with the multiplicative version. In what follows, we will present only results which correspond to Lemma 3.2, Theorem 3.1, 3.2, and Corollary 3.1 and omit the proofs since they are almost identical with those presented above. Notice that the other results are also valid in this setting, thought we will not formulate explicitly.

Let us define the linear operator $U_{\psi}$ by
\[U_{\psi} : L^{\infty}(\mathbb{R}_+^{\times}) \longrightarrow L^{\infty}(\mathbb{R}_+^{\times}), \quad  (U_{\psi}f)(x) = \int_1^x f(t)\psi(x/t)\frac{dt}{t}. \]

\begin{lem}
For every $\psi \in L^1(\mathbb{R}_+^{\times})$, $f \in L^{\infty}(\mathbb{R}_+^{\times})$ and $\omega \in \Omega^*$, it holds that
\[e^{\omega}\mathchar`-\lim_x \int_1^x f(t)\psi(x/t)\frac{dt}{t} = \int_1^{\infty} f_{\omega}^{\times}(t^{-1})\psi(t)\frac{dt}{t}. \]
\end{lem}

\begin{thm}
For every $\psi \in L^1(\mathbb{R}_+^{\times})$, $f \in L^{\infty}(\mathbb{R}_+^{\times})$ and $\omega \in \Omega^*$, we have
\[(U_{\psi}f)_{\omega}^{\times}(x) = (f_{\omega}^{\times} * \tilde{\psi})(x) \]
for every $x \in \mathbb{R}^{\times}$, where $\tilde{\psi}$ is the function in $L^1(\mathbb{R}^{\times})$ such that
\[\tilde{\psi}(x) = \begin{cases}
            \psi(x) & \text{if $x \ge 1$}, \\
            0       & \text{if $0 < x <1$}.
         \end{cases}     
\]
\end{thm}

\begin{thm}
For every $\psi \in L^1(\mathbb{R}_+^{\times})$ with $\int_1^{\infty} \psi(x)\frac{dx}{x} = 1$, the summability method $M_{\psi}$ is equivalent to the summability method $P$ if and only if the closed ideal $J$ generated by $\tilde{\psi}$ in $L^1(\mathbb{R}^{\times})$ is the whole space $L^1(\mathbb{R}^{\times})$.
\end{thm}

\begin{cor}
For every $\psi \in L^1(\mathbb{R}_+^{\times})$ with $\int_1^{\infty} \psi(x)\frac{dx}{x} = 1$, the summability method $M_{\psi}$ is equivalent to the summability method $P$ if and only if $Z(\tilde{\psi}) = \emptyset$.
\end{cor}

Now let us consider some examples. Notice that the characters of $\mathbb{R}^{\times}$ are $t^{ix} = e^{ix\log t} (x \in \mathbb{R})$.

First we put $\psi_r(x) = rx^{-r} (r > 0)$. Let us denote $M_{\psi_r}$ by $M_r$, namely,
\[M_r(f) = \lim_{x \to \infty} \int_1^x f(t) r\left(\frac{x}{t}\right)^{-r}\frac{dt}{t} = \lim_{x \to \infty} \frac{r}{x^r} \int_1^x f(t)t^{r-1}dt, \]
where $f \in L^{\infty}(\mathbb{R}_+^{\times})$. In particular, for $r=1$, we have Ces\`{a}ro method $M = M_1$.
Now the direct computation shows that
\[\hat{\psi}_r(x) = \int_{\mathbb{R}^{\times}} rt^{-r}t^{-ix}\frac{dt}{t} = \frac{r}{r+xi} \quad (r > 0).  \]
Hence $\hat{\psi}_r(x)$ has no zero point in $\mathbb{R}$ for every $r > 0$ and hence each $M_r$ is equivalent to $P$. Thus we get the following result immediately.

\begin{thm}
For any $f \in L^{\infty}(\mathbb{R}_+^{\times})$, 
\[M(f) = \lim_{x \to \infty} \frac{1}{x} \int_1^x f(t)dt = \alpha \]
holds if and only if 
\[M_r(f) = \lim_{x \to \infty} \frac{r}{x^r} \int_1^x f(t)t^{r-1}dt = \alpha \]
holds for some $r > 0$.
\end{thm}

Remark that this type of summability methods were considered in the discrete setting by several authors, for example, [2] and [5]. In particular, the following result was stated in [5] as a consequence of the result in [1].
\begin{thm}
For any $f \in l_{\infty}$, 
\[\lim_{n \to \infty} \frac{1}{n} \sum_{i=1}^n f(i) = \alpha \]
holds if and only if
\[\lim_{n \to \infty} \frac{r}{n^r} \sum_{i=1}^n f(i)n^{r-1} = \alpha \]
holds for some $r > 0$.
\end{thm}

Notice that this result can be derived from Theorem 3.7 by using the natural embedding $V : l_{\infty} \rightarrow L^{\infty}(\mathbb{R}_+^{\times})$ defined by $(Vf)(x) = f([x])$.

Next we consider iterations of $M$. Let us $(Uf)(x) = \frac{1}{x} \int_1^x f(t)dt$. We define a summability method $H_2$ by
\[H_2(f) = \lim_{x \to \infty} \frac{1}{x} \int_1^x \left(\frac{1}{t} \int_1^t f(s)ds\right)dt = \lim_{x \to \infty} (U^2f)(x). \]
In this way, we define inductively the sequence of summability methods $H_1 = M, H_2, \ldots, H_k, \ldots$, where
\[H_k(f) = \lim_{x \to \infty} (U^kf)(x). \]
These summability methods were considered by H\"{o}lder in the discrete setting (see [3]). Then by Corollary 3.4, we have the following result.

\begin{thm}
Every $H_k, k \ge 1$, is  equivalent to $H_1 = M$.
\end{thm}

\section{Some other types of summability methods}
For any $\phi \in L^1(\mathbb{R}_+)$ with $\int_0^{\infty}\phi(x)dx = 1$, we define a summability method $S^*_{\phi}$ as follows:
\[S^*_{\phi}(f) = \lim_{x \to \infty} \int_x^{\infty} f(t)\phi(t-x)dt, \]
where $f(x) \in L^{\infty}(\mathbb{R}_+)$. We define the linear operator $U_{\phi}^*$ by 
\[U_{\phi}^* : L^{\infty}(\mathbb{R}_+) \longrightarrow L^{\infty}(\mathbb{R}_+), \quad  (U_{\phi}^*f)(x) = \int_x^{\infty} f(t)\phi(t-x)dt. \]

It is easy to show that $(U_{\phi}^*f)(x) \in C_{ub}(\mathbb{R}_+)$. Then we can prove the following result in the same way as the previous section.
\begin{lem}
For every $\phi \in L^1(\mathbb{R}_+)$, $f \in L^{\infty}(\mathbb{R}_+)$ and $\omega \in \Omega^*$, it holds that
\[\omega\mathchar`-\lim_x \int_x^{\infty} f(t)\phi(t-x)dt = \int_0^{\infty} f_{\omega}(t)\phi(t)dt. \]
\end{lem}

\begin{thm}
For every $\phi \in L^1(\mathbb{R}_+)$, $f \in L^{\infty}(\mathbb{R}_+)$ and $\omega \in \Omega^*$, we have
\[(U_{\phi}^*f)_{\omega}(x) = \int_{-\infty}^{\infty} f_{\omega}(x+t)\tilde{\phi}(t)dt = (f_{\omega} * \phi^*)(-x)  \]
for every $x \in \mathbb{R}$, where $\phi^*$ is a function in $L^1(\mathbb{R})$ such that
\[\phi^*(x) = \begin{cases}
            0 & \text{if $x \ge 0$}, \\
            \phi(-x)  & \text{if $x < 0$}.
         \end{cases}     
\]
\end{thm}

\begin{thm}
For every $\phi \in L^1(\mathbb{R}_+)$, the summability method $S_{\phi}^*$ is equivalent to the summability method $K$ if and only if the closed ideal $I^*$ generated by $\phi^*$ in $L^1(\mathbb{R})$ is the whole space $L^1(\mathbb{R})$.
\end{thm}

Notice that $\phi^*(x) = \tilde{\phi}(-x)$ and
\[\hat{\phi^*}(x) = \int_{-\infty}^{\infty} \phi^*(t)e^{-itx}dt = \int_{-\infty}^{\infty} \tilde{\phi}(-t)e^{-itx}dt = \int_{-\infty}^{\infty} \tilde{\phi}(t)e^{itx}dt = \hat{\tilde{\phi}}(-x). \]
This observation implies that $I^* = L^1(\mathbb{R})$ if and only if $I = L^1(\mathbb{R})$. Then we have the following result.

\begin{thm}
For every $\phi \in L^1(\mathbb{R}_+)$ with $\int_0^{\infty} \phi(x)dx = 1$ and $Z(\tilde{\phi}) = \emptyset$, the summability methods $S_{\phi}$ and $S_{\phi}^*$ are equivalent. Namely, for every $f(x) \in L^{\infty}(\mathbb{R}_+)$, it holds that
\[\lim_{x \to \infty} \int_0^x f(t)\phi(x-t)dt = \alpha \]
if and only if
\[\lim_{x \to \infty} \int_x^{\infty} f(t)\phi(t-x)dt = \alpha. \]
\end{thm}

In the same way, for any $\psi \in L^1(\mathbb{R}_+^{\times})$ with $\int_1^{\infty} \psi(x)\frac{dx}{x} = 1$, we define a summability method $M_{\psi}^*$ as follows:
\[M_{\psi}^*(f) = \lim_{x \to \infty} \int_x^{\infty} f(t)\psi(t/x)\frac{dt}{t}, \]
where $f(x) \in L^{\infty}(\mathbb{R}_+^{\times})$, and then we have the following analogous result.

\begin{thm}
For every $\psi \in L^1(\mathbb{R}_+^{\times})$ with $\int_1^{\infty} \psi(x)\frac{dx}{x} = 1$ and $Z(\tilde{\psi}) = \emptyset$, the summability methods $M_{\psi}$ and $M_{\psi}^*$ are equivalent. Namely, for each $f(x) \in L^{\infty}(\mathbb{R}_+^{\times})$, it holds that
\[\lim_{x \to \infty} \int_1^x f(t)\psi(x/t)\frac{dt}{t} = \alpha \]
if and only if
\[\lim_{x \to \infty} \int_x^{\infty} f(t)\psi(t/x)\frac{dt}{t} = \alpha. \]
\end{thm}

Combining Theorem 3.7 with Theorem 4.4, we obtain the following corollary. 

\begin{cor}
For any $f(x) \in L^{\infty}(\mathbb{R}_+^{\times})$ and $r > 0$, we have
\[\lim_{x \to \infty} \frac{r}{x^r} \int_1^x f(t)t^{r-1}dt = \alpha \]
if and only if
\[\lim_{x \to \infty} rx^r\int_x^{\infty} f(t)\frac{dt}{t^{r+1}} = \alpha. \]
\end{cor}

\end{document}